\documentclass[%
 reprint,
 amsmath,amssymb,
 aps,
 pre,
showkeys
]{revtex4-2}

\usepackage{graphicx}% Include figure files
\usepackage{dcolumn}% Align table columns on decimal point
\usepackage{bm}% bold math
\usepackage{float}
\usepackage{cancel}
\usepackage{xcolor}

\usepackage{amsmath}
\newcommand\numberthis{\addtocounter{equation}{1}\tag{\theequation}}
\begin{document}
\preprint{APS/123-QED}

\title{Sloppy model analysis identifies bifurcation parameters without normal form analysis}

\author{Christian N. K. Anderson}
 \email{bifurcate@byu.edu}
\author{Mark K. Transtrum}
\email{mktranstrum@byu.edu}
 \affiliation{Department of Physics and Astronomy\\ Brigham Young University\\
 Provo, UT 84604.}

\date{\today}% It is always \today, today

\begin{abstract} Bifurcation phenomena are common in multi-dimensional multi-parameter dynamical systems.
Normal form theory suggests that bifurcations are driven by relatively few combinations of parameters.
Models of complex systems, however, rarely appear in normal form, and bifurcations are controlled by nonlinear combinations of the bare parameters of differential equations.
Discovering reparameterizations to transform complex equations into a normal form is often very difficult, and the reparameterization may not even exist in a closed-form.
Here, we show that information geometry and sloppy model analysis using the Fisher information matrix can be used to identify the combination of parameters that control bifurcations.
By considering observations on increasingly long time scales, we find those parameters that rapidly characterize the system's topological inhomogeneities, whether the system is in normal form or not.
We anticipate that this novel analytical method, which we call time-widening information geometry (TWIG), will be useful in applied network analysis.
\begin{description}
\item[Software] An implementation and example of TWIG analysis are available at https://github.com/oceanchaos/TWIG/
\end{description}
\end{abstract}


\maketitle

\section{\label{sec:level1}Introduction}

This paper provides a method for extracting bifurcation parameters from a set of dynamic equations by combining information geometry and bifurcation theory. Both are useful for modeling multi-parameter systems and systems with multiple regimes of behavior respectively, but together they provide methods for data-driven analysis of a wide array of natural phenomena. By creating an explicit connection between the information in the signal (model output) and the model parameters, we identify the combinations of parameters responsible for topological change in the dynamics, the codimension of the bifurcation, and the time scale necessary to resolve this information. The information further provides the directions normal to the separatrix, which divides behavioral regimes of the system.

%Research into bifurcation phenomena continues to grow. Initially neglected after Poincaré's investigations in 1890, it experienced a renaissance in the 1960s~\cite{Holmes1990} and is now a bustling field of research with over 5,000 articles using the keyword in 2020 alone.\footnote{According to ISI Web of Science, search performed May 2021}

Traditionally, when confronted with a high-dimensional, multi-parameter system of dynamic equations, bifurcation analysis proceeds by attempting to simplify the system to a manageable size.
Center Manifold Reduction exploits the Hartman-Grobman theorem~\cite{Hartman1982} to create a lower-dimensional linear map in the region of a critical point that is locally accurate and is a rapid way to determine the system stability.
Shoshitaivishili extended this method to non-hyperbolic equilibria, creating a container for critical modes to straighten out non-linear terms and, ideally, drop some of them~\cite{Crawford1991}. Such methods have been used to describe phenomena as diverse as neural network optimization and foraging decisions in monkeys~\cite{Brown2005, Feng2009}.

A related approach is the method of Poincaré-Birkhoff normal forms. 
It uses appropriately centered manifolds to analyze which nonlinear terms are essential and must remain even under optimal coordinate transformations. 
Such transformations are useful, because the reduced normal-form equations typically have greater symmetry than the initial problem, a property that can be exploited by many analytical tools. 
Though powerful, ``in practice lengthy calculations may be necessary to extract the relevant normal-form coefficients from the initial equations ~\cite{[Quote from p. 1021 of ]Crawford1991}.''
Even if such coefficients can be found, neither their interrelationship nor their relative sensitivities are always apparent. 
It is often the case that some parameters differ by many orders of magnitude in their effect on long-term dynamics, and a method that doesn't distinguish among them is sub-optimal for most applications.

The method of Lyupanov exponents is an admirably general tool for analyzing the global stability of a system. 
Unfortunately, it provides little information about which specific parameter combinations lead to system (in)stability.
For the purposes of bifurcation analysis, it is therefore sometimes paired with sensitivity analyses based on the global sensitivity metrics of Sobol'~\cite{Sobol}.
These measures, along with useful extensions such as FAST (Fourier amplitude sensitivity test) and Importance Measures~\cite{Bettonvil1997,Homma1996,Gul2015}, are able to determine exactly how much of a model's variability is due to each of its parameters.
While this often works in practice, there are two potential pitfalls in this approach.
First, it assumes that the parameters responsible for variability are also responsible for instability, which is not always the case.
Second, if the bifurcation is caused by combinations of many parameters (as frequently happens), then variability will often be high across all these parameters even though the bifurcation itself has a low codimension. 
In other words, a low-dimensional bifurcation surface generally cuts diagonally across parameter space unless appropriately reparameterized.
Once such a transformation is applied and the system is reduced to a normal form (see Sec.~\ref{sec:normform}), then the codimension should be apparent, but finding that reparameterization is still likely to be cumbersome, if not impossible, in closed form. 
Just one such transformation can require several papers, as in the case of high-dimensional diffusion-activated processes from Kramers, through Langer, and finally to one dimension, derived using iterations of singular value decomposition by Berezhkovskii~\cite{Berezhkovskii2005}.

A third, independent line of analysis comes from Renormalization Group (RG) methods, which are usually applied to study universal power laws near critical points.
Feigenbaum~\cite{Feigenbaum1978} was the first to note such universalities in bifurcations of the discrete period-doubling type, a result that he and others extended until it included all major bifurcation types~\cite{Feigenbaum1979,Widom1982,Hu1982,Hu1995,Hathcock2021}.
Working from the other direction, scientists investigating critical phenomena with RG (e.g., many behaviors of quantum chromodynamics) have discovered bifurcations, and used the tools of one to analyze the other~\cite{Chang1985, raju2019normal}.
A remarkable study found deep equivalence between RG transformations and normal form theory, showing that the difficult transformation of an ODE system into a normal form could often be accomplished to at least second order by applying three RG transforms~\cite{DeVille2008}.

More broadly, universal scaling laws and RG analysis of critical points is often associated with emergence and the systematic irrelevance of many degrees of freedom.
Recent work has extended these ideas to a broader class of systems known as ``sloppy models''~\cite{Transtrum2015,Transtrum2010,Transtrum2016,Machta2013,Brown2004a}.
The moniker ``sloppy'' is meant to convey that these systems have a few combinations of parameters that are many orders of magnitude more influential than other parameter combinations. 
More precisely, one unit change in a ``stiff'' parameter direction has as much influence as a million or more unit change in a different ``sloppy'' direction. 
Sloppy model analysis relies heavily on the techniques of information geometry~\cite{Transtrum2010,Transtrum2011,quinn2022information} and in this paper we use the terms interchangeably.
These techniques have motivated  novel reduction algorithms by removing unimportant, sloppy parameters~\cite{Transtrum2014, mattingly2018maximizing, quinn2022information}.
Recent work~\cite{Raju2018} demonstrates that as coarse-graining of RG models proceeds, the flow causes information of ``relevant'' parameter combinations to be maintained while ``irrelevant'' parameters are compressed and become sloppy.
These ideas share a common goal with bifurcations analysis in which many diverse systems are collected into a few universal, normal forms.
This paper closes the loop, showing how information geometry applies directly to bifurcation analysis without passing through the ``middleman'' of renormalization group theory. 
The usefulness of such an analysis, which we call Time Widening Information Geometry (TWIG), also circumvents the need for the other types of analyses described above.

In this work, we demonstrate similar notions of ``relevant'' and ``irrelevant'' parameters near a bifurcation using the formalism of information geometry and sloppy models.
The intuition behind this approach is as follows.
Topological inhomogeneities in the flow field produce trajectories containing different information on either side of a bifurcation.
For example, on one side of a Hopf bifurcation, all trajectories collect into a central fixed point, while they flow into an orbit (limit cycle) on the other side.
TWIG works by measuring the information content in these trajectories at increasingly long time scales and identifying those combinations of parameters to which the trajectory is most sensitive. 
At long time scales, these are the parameters responsible for the bifurcation, while parameters that cause only local variability have less impact.

Information geometry can be applied to complex systems from many disciplines--but especially systems biology--to iteratively ``reverse engineer'' optimal statistical models by removing parameters whose value has little influence on the macroscopic behavior of the system~\cite{Transtrum2014, Transtrum2015, Transtrum2016, Jeong2018}. However, it was recently appreciated that such reverse engineering can be done even if the underlying system bifurcates into qualitatively different behaviors, because the information geometry of parameters participating in the bifurcation show a characteristic ``sand dune'' shape when crossing from one behavioral state to another~\cite{Roesch2019}. These results imply that if the functional form of the system is known, it should be even easier to determine bifurcation parameters than if the system's equations need to be inferred.

This paper is organized as follows: In Section~\ref{sec:bkg}, we provide background information on bifurcations and information geometry generally, and, specifically, how we conceptualize them for the purposes of applying the latter to the analysis of the former. In Section~\ref{sec:normform}, we show how an IG analysis of the normal form bifurcations rapidly provides insight into the structure of bifurcations simple enough to be understood by other methods. Section~\ref{sec:abnorm} shows how this analysis extends to more difficult bifurcations, the implications of which are summarized in Section~\ref{sec:conc}.

\section{\label{sec:bkg}Background and Problem Formulation}

\subsection{\label{sec:bifintro}Bifurcations}

Bifurcations frequently arise in the analysis of dynamical systems, where one typically characterizes the flow field with special attention to any fixed points or stable oscillations~\cite{Strogatz2015}. Consider a generalized system of $n$ coupled dynamic equations, where each equation is of the form $\dot{\mathbf y} = f(\mathbf y; \mathbf \theta)$, where $\mathbf \theta$ is a vector of $m$ parameters. Small changes to any of the $\theta_i$ values typically result in correspondingly small changes to the $n$-dimensional vector field, such as small changes to the position of a fixed point or radius of a limit cycle. Such deformations are topologically equivalent (meaning the number and properties of the attractors / repellers in the field do not change) and  homeomorphic (continuous with a continuous inverse). However, there may be critical parameter values where a small change causes new fixed points to emerge from old ones, or two fixed points to approach and be mutually annihilated, or limit cycles to be broken. Such nonhomeomorphic transformations are generically called bifurcation no matter their exact form.

Several types of simple bifurcations have been identified and reduced to their simplest possible mathematical expression. These are the so-called ``normal forms" and are enumerated in the section below. These forms are convenient starting points for analysis, since they have clearly defined rate parameters that are unambigiously responsible for causing topological inhomogeneities. However, even elegant mathematical descriptions of real-world dynamical systems rarely conform exactly to one of the normal forms. 

Bifurcation parameters for physical models often do not align with the bare parameters. In the classic example of boiling liquid, the bifurcation parameter is some combination of temperature, pressure, salinity, and others. In general, a reparameterization to a single, unambiguous bifurcation parameter may be possible in principle, but often requires either substantial additional physical insight, or mathematical sophistication, or both. Some researchers have even recommended building an analogous physical circuit as the fastest method to detect the bifurcation~\cite{JimenezRamirez2021}. Complex models can have hundreds of coupled dynamic equations with thousands of parameters (e.g., models of sophisticated mobile phone circuit boards~\cite{Kleijnen2006}, or complex protein networks~\cite{Waldherr2007}). How can we determine which parameter (or more likely, combination of parameters) is responsible for the bifurcation in such cases? 

\begin{figure*}
\includegraphics{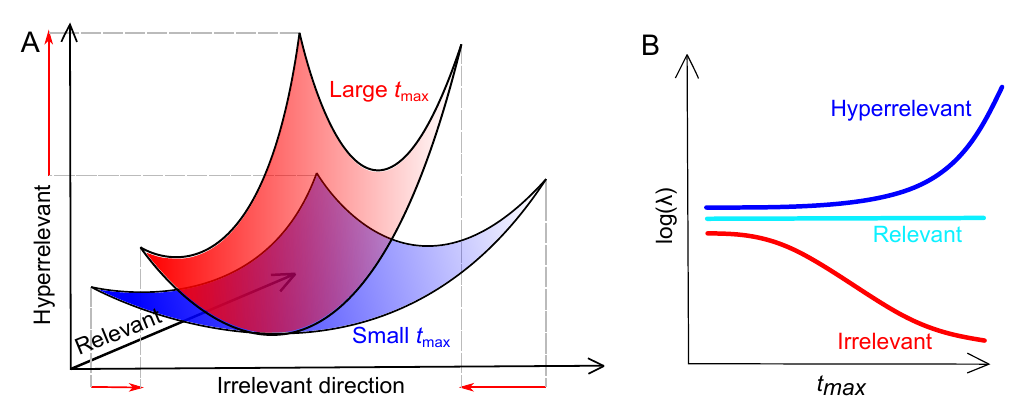}
\caption{(A) The model manifold in data space represents all values that can be reached. The axes represent directions that are distorted in characteristic ways as $t_{max}$ increases. They can be contracted (irrelevant), expanded (hyperrelevant) or unchanged (relevant). (B) Relevance can be quantified by observing the eigenvalues of the Fisher information matrix as $t_{max}$ is increased. Eigenvalues that do not change at longer time scales retain their relevance, while those that increase or decrease become either more or less relevant.}
\label{fig:eigfig}
\end{figure*}

\subsection{\label{infointro}Information Geometry}

The fundamental object of information geometry is the Fisher information matrix (FIM or $\mathcal I$), which quantifies the information that the observations $\mathbf y$ contain about the parameters $\mathbf \theta$ of a dynamical system. Here we introduce the FIM for dynamical systems.

Consider a system of ordinary differential equations where the parameters are tuned to be exactly at their critical values, i.e., the system is at (one of) its bifurcation point(s). The system is allowed to evolve, and the trajectory of one of its equations $y_j$ is sampled at several time points $y_{j}(t_i)$, where $t_i = t_0 + \frac{i}{n} t_{max}$. 
%The solution to the differential equation $\dot{y}=-\theta_2 e^{-\theta_2 t}+\theta_3 e^{\theta_3 t}$ is:
To help visualize this process, let us imagine a one-dimensional system
\begin{equation}
    y(t)=\theta_1+e^{-\theta_2 t}+e^{\theta_3 t}
    \label{eq:toymodel}
\end{equation} sampled at $t=\{1,2,3\}$ to create a vector of three observations $\mathbf{y}=\{y(t_1),y(t_2),y(t_3)\}$ which we plot in $\mathbb{R}^3$, i.e., data space. If $\theta_3>0$, then there is no equilibrium; if $\theta_3=0$ and $\theta_2>0$ then the equilibrium is at $\theta_1+1$ or $\theta_1+2$ if $\theta_2=0$. As the parameters of $\mathbf{\theta}$ change, the position of $\mathbf{y}$ will also change, but except for extreme values of $\theta_i$, it cannot reach all possible values in $\mathbb{R}^3$. The space filled by values of $\mathbf{y}$ that can be reached for a given range of parameter values defines the model manifold. A schematic of such a manifold is drawn in Fig.~\ref{fig:eigfig}A. 

The Fisher information is most commonly defined in probabilistic terms as the expected Hessian matrix of the log-likelihood: 
\begin{equation}
  \mathcal{I}=-E\left[\frac{\partial^2}{\partial \theta^2} \log \mathcal{L}(\theta|d)\right]
    \label{eq:FIM}
\end{equation}
where $\theta$ is a vector of parameters, and $d$ is the data.
For deterministic systems such as we consider here, it is standard practice to assume that measurements are obscured by additive Gaussian noise,
\begin{equation}
  \label{eq:noise}
  d_i = y(t_i) + \zeta
\end{equation}
where $y(t_i)$ is the (deterministic) output of the model at time $t_i$ and $\zeta$ is standard normal random variable $\zeta \sim \mathcal{N}(0, 1)$.
This assumption defines a probability distribution to which Eq.~\ref{eq:FIM} can be applied~\cite{Transtrum2011}.
Because this construction is so common in information theory, it is often referred to as the sensitivity Fisher information matrix or sFIM~\cite{ Brouwer2018} for reasons that will soon be apparent.
% However, in this context it isn't immediately clear what the probability function underlying the likelihood would be. 
% Additionally, in practice, the second derivatives are time-consuming and difficult to calculate. 
In general $\mathcal{I}$ can be expressed in terms of the first derivatives only
\begin{align}
  \label{eq:FIMrelations}
  \mathcal{I} & =-E\left[\frac{\partial^2}{\partial \theta^2} \log \mathcal{L}(\theta|d)\right] \\
              & = E\left[\frac{\partial}{\partial \theta} \log \mathcal{L}(\theta|d) \frac{\partial}{\partial \theta} \log \mathcal{L}(\theta|d)\right]
\end{align}
Using the second form, one can show that sFIM becomes
\begin{align}
  \label{eq:JTJ}
  \mathcal{I}_{i,j} & = \sum_{k = 1}^M J_{k,i} J_{k,j} = (J^T J)_{i,j}
\end{align}
where we have introduced the Jacobian or sensitivity matrix $J_{k,j}=\frac{\partial y_k}{\partial \theta_j}$ whose entries denote the sensitivity of prediction $k$ to changes in parameter $j$.
In Eq.~\eqref{eq:JTJ}, $M$ denotes the number of observations.

The entries of the FIM indicate the sensitivity of the model's trajectory to changes in each pair of parameters. A high score indicates that a parameter pair has a strong influence on model dynamics, while a small score indicates a “sloppy” direction (parameter values can change a great deal without much changing $\mathbf y$). 
The curvature of the likelihood function converts distances in parameter space to distances on the manifold in data space, making the FIM a Riemannian metric tensor on the model manifold in data space.
It is important to note that the physical units of parameters can strongly affect the values of the FIM.
For this reason, it is common to perform dimensional analysis before sloppy model analysis as we do throughout this study.

In general the curvature of the likelihood surface does not align with the bare parameters.
Rather, the characterization of the model's sloppiness aligns with the eigenvectors of $\mathcal{I}$. Eigenvalues of the FIM are related to the singular value decomposition of $J=U\Sigma V^T$:
\begin{equation}
    \begin{aligned}
        \mathcal{I}&=V\Sigma^2V^T.
    \end{aligned}
    \label{eq:SVD}
  \end{equation}
where $U$ and $V$ are matrices of the left and right singular vectors of $J$, and $\Sigma$ is the diagonal matrix of its singular values. This implies that the right singular vectors of the Jacobian $V$ are also the eigenvectors of the FIM.
The eigenvectors of $\mathcal{I}$ ``orient'' the parameter-space into the parameter combinations most relevant for changing the model's behavior.

Imagine now that we coarsen the sampling rate by changing $t_{max}$.
In our simple example, increasing $t_{max}$ from 3 to 6 means the model $\mathbf{y}$ is sampled at $t=\{2,4,6\}$.
This procedure stretches the manifold in some directions and compresses it in others.
This distortion is measured by an increase or decrease in the eigenvalues of $\mathcal{I}$, respectively.
Compression of the manifold (i.e., decreasing eigenvalue) with increasing $t_{max}$ indicates that the combination of parameters is less important to the long-term dynamics.
We call the corresponding eigendirection ``irrelevant".
Similarly, if the manifold stretches (i.e., increasing eigenvalue), we call the corresponding direction ``hyperrelevant".
Directions that are neither compressed nor stretched are called ``relevant" directions (Fig.~\ref{fig:eigfig}B).
Returning to the example in Eq.~\ref{eq:toymodel}, $\theta_1$ is relevant since its effect on the model's output is unchanged with observation time.
In contrast, $\theta_2$ is irrelevant since the exact rate of the decay matters less as time scales become very large, and $\theta_3$ is hyperrelevant since small changes have large effects at large $t$. Note that $\theta_2$ and $\theta_3$ are functionally interchangeable if either is negative.

This procedure is similar to coarse-graining under RG flow described in reference~\cite{Raju2018} and is used to generate their Fig.~1. In our case, however, because we are coarsening the sampling rate, the total observation time increases and includes new information, i.e., observations at later times.  As such, it is not a true coarse-graining and introduces the possibility of hyperrelevant directions, i.e., directions that become increasingly important such as $\theta_3$.
We will see that hyperrelevant directions are associated with the stability or instability of the equilibrium.

This method is also somewhat analogous to studies that use Sobol' sensitivity analysis to track importance at different time scales, either bare parameters or eigenvalue combinations. Such methods are excellent at providing estimates of model variability at a given point in parameter space, and have noted both increasing and decreasing importance for model parameters of biophysical systems~\cite{Song2013,Alexanderian2020}. Critics note that these methods are computationally expensive, even when implementing Morris acceleration~\cite{Sumner2012}, and the implications for bifurcation analysis are not immediately obvious.

In addition to characterizing bifurcations, TWIG analysis reveals two other features of bifurcating systems. First, there can be parameters (or combination of parameters) that move the location of a fixed point without causing a bifurcation. 
Such parameter combinations appear as ``relevant'' eigendirections, as the new equilibrium appears in long-time observations. 
These parameters need to be removed in order to correctly identify the codimension of the bifurcation.
We do this by solving for the location of the fixed point with a numerical root-finding algorithm and subtracting it from the trajectory at every point. This effectively translates the fixed point to the origin and is analogous to the recentering step of Center Manifold Analysis. 
For limit cycle trajectories, we recenter by subtracting off the (unstable) fixed point that must exist within the cycle (according to the Poincar\'{e}-Bendixson theorem~\cite{Bendixson1901}).

\begin{figure}[b]
\includegraphics[width=.98\linewidth]{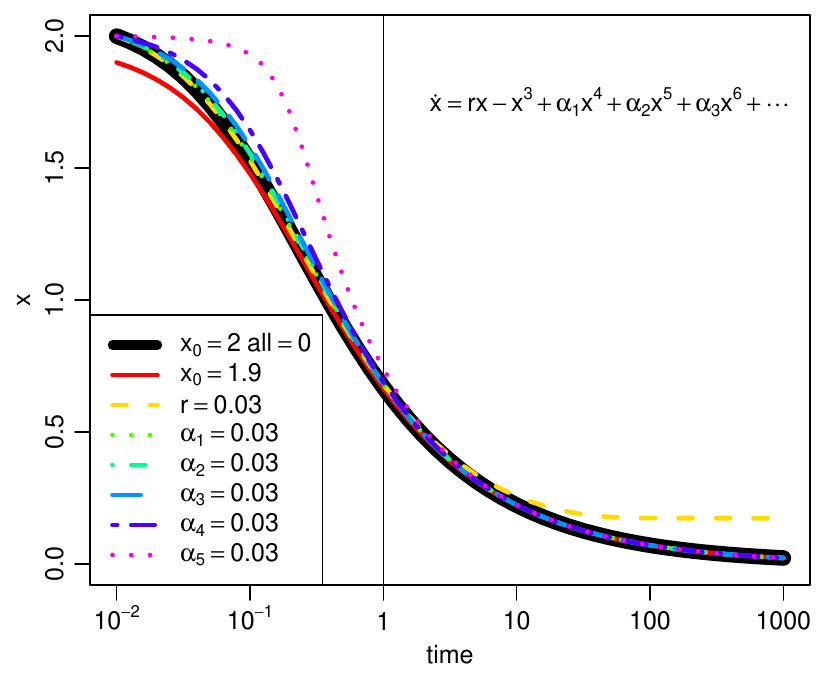}
\caption{\label{fig:pitchtraj}The trajectory of a supercritical pitchfork at the bifurcation point (heavy black line), and slightly perturbed from it (thin colored lines). At short time scales, $y_0$ (thin red) and high-order parameters (long-dashed blue) appear relevant.
But, as the dynamics progress, $r$ (short-dashed yellow) emerges as the only parameter that changes the long-term equilibrium point. 
This change from relevant to not (and vice versa) occurs at $t \approx 10$, and is reflected in the arch shape and changing colors of Fig.~\ref{fig:pitcheigs}}.
\end{figure}

The second feature arises in such oscillating systems.  Parameters that change the phase or frequency of oscillation without destroying the equilibrium itself appear as hyperrelevant as the accumulating phase difference becomes increasingly important at late times.  Previous research has shown that such systems frequently cause problems in an information geometry framework by introducing ``ripples" into the likelihood surface of Eq.~\ref{eq:FIM}. The solution is to perform a coordinate transformation so the period itself becomes a parameter. In one formulation of the FIM, this causes the manifold to ``unwind'', creating a smooth likelihood surface~\cite{Francis2019} and thereby eliminating a misleading eigendirection. 
%Here, transforming the Hopf into polar coordinates clarified its analysis.

Four important pieces of information come from this Time Widening Information Geometry (TWIG) analysis. First, the number of hyperrelevant and relevant directions corresponds to the codimension of the bifurcation system. Second, the square of each element of the eigenvector matrix $V_{ij}$ indicates the participation factor of each bare parameter $\theta_i$ in eigenvector $j$. This last fact follows because the participation factor $p_{ij} \equiv U_{ij} V_{ij} =V_{ij}^2$ as can be seen by combining the definition of a participation factor~\cite{Perez-Arriaga1982, Garofalo2002} with Eq.~\ref{eq:SVD} above which shows that the left and right singular vectors of the FIM are identical. Third, the eigendirections themselves will change as $t_{max}$ increases and parameters that influence the short-term dynamics lose their salience at long time scales. If initial conditions are included as parameters, their loss of relevance is a strong indicator that the system has been simulated ``long enough" to capture equilibrium behavior. This is not a trivial concern in practice, where long numerical simulations are always fighting the accumulation of computer round-off error. Finally, at equilibrium the relevant eigendirections point along the (potentially) high-dimensional separatrix surface, and so the bifurcation can be mapped through all parameter space. 

Note that this procedure works no matter the number of dynamical variables involved in the differential equation system.
However, it presupposes that the model can be simulated on at least one side of the bifurcation to arbitrarily long times, \emph{i.e.} it analyzes stable dynamics on the threshold of instability.
A bifurcation that switches between two different forms of instability will not be easily detectable with this method, as trajectories will diverge on both sides of the bifurcation.
In the next section, we demonstrate how this procedure works for all common normal forms of bifurcations.

\section{\label{sec:normform}Normal-form Bifurcations}

Local bifurcations can be described mathematically in a potentially infinite number of ways, but nearly all of them can be reparameterized, at least locally, to one of five kinds of normal forms. These are:
\begin{itemize}
    \item Saddle-node: $\dot{x} =r +x^2$, where one stable and one unstable fixed point emerge from a previously uninterrupted flow at a critical value $r_{crit}$
    \item Transcritical: $\dot{x} = r x-x^2$, where a stable and an unstable fixed point exist everywhere but at the bifurcation, and swap stability at the critical value
    \item Supercritical Pitchfork: $\dot{x} = rx-x^3$, where symmetric stable fixed points emerge from a single fixed point, which itself becomes unstable
    \item Subcritical Pitchfork: $\dot{x} = rx+x^3$, symmetric unstable fixed points emerge from an unstable fixed point, which itself becomes stable
    \item Hopf: a stable limit cycle emerges from what had previously been a stable point attractor. Depending on the coordinate system, the normal form is $\dot{z}=z(a+b|z|^2)$ (complex), $\dot{x} = -y + x(\mu-r^2);~\dot{y} = x + y(\mu-r^2)$ (Cartesian), or $\dot{r} = r(\mu-r^2); ~\dot\theta = -1$ (Polar). 
\end{itemize}
\begin{figure}
\includegraphics[width=.98\linewidth]{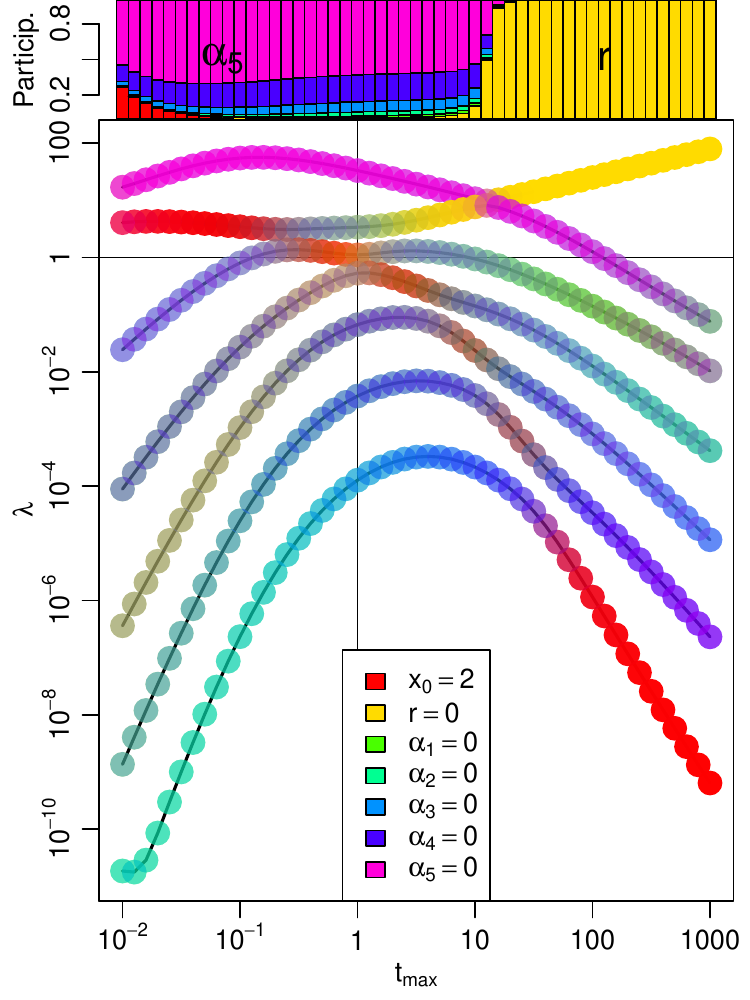}
\caption{\label{fig:pitcheigs}The ``rainbow diagram'' of the system from Fig.~\ref{fig:pitchtraj}, showing the FIM's eigenanalysis at each $t_{max}$. The top panel represents the participation of each parameter in the first eigenvector ($V_{i,1}^2$ in Eq.~\ref{eq:SVD}). The leading eigenvector changes from $\alpha_5$-dominant (purple) to $r$-dominant (yellow) at $t_{max} \approx 10$, \emph{i.e.}, just where the $\alpha_5$ trajectory is replaced by the $r$ trajectory as most divergent in Fig.~\ref{fig:pitchtraj}. 
The large panel below shows all seven eigenvalues ($\lambda_i = \sqrt{\Sigma^2_{ii}}$ in Eq.~\ref{eq:SVD}) at each $t_{max}$, colored as the weighted average RGB of each parameter's participation factor. Thus, the top line, corresponding to the largest eigenvalue in the top panel, starts mostly purple ($\alpha_5$), but turns yellow as $r$ dominates the first eigenvalue at larger $t_{max}$ values. For all parameters, a small change to parameter values influences trajectory at short time scales (the rising limb) but, with the exception of $r$, \emph{not} at long time scales (the descending limb). The red color in the bottom-right indicates that the initial value $x_0$ eventually becomes the least relevant parameter in the model. Pure colors indicate an eigenvector pointing along a parameter axis, while mixed colors like browns and greys indicate that many parameters participate in the eigenvector.}
\end{figure}
A method able to detect bifurcation parameters for these types of bifurcations will detect the overwhelming majority of bifurcations we are likely to encounter. 
The Fisher information as a function of $t_{max}$ for each bifurcation type has a closed-form solution, which complements and validates the numerical results that we present here (see Appendix~\ref{apx:SN}).
In each case, the sensitivity with respect to the bifurcation parameter, $r$, dominates the long-term dynamics of the system in the neighborhood of the bifurcation, no matter how many other higher order parameters are added to the normal form.

For example, a supercritical pitchfork of the form $\dot{x} = r x-x^3+\alpha_1 x^4+\alpha_2 x^5\dots$ experiences a bifurcation when $r=\alpha_i=0$.
At short time scales (e.g., where $t_{max}<1$), the system's trajectory is strongly influenced by changes to its initial value $x_0$ and the higher order $\alpha$ terms (for $x_0>1$).
However, later dynamics show that changes to the $\alpha_i$'s (and $x_0$) barely affect the trajectory of approach to equilibrium at 0, while small modifications to $r$ move the equilibrium itself (Fig.~\ref{fig:pitchtraj}).
An eigenanalysis of the FIM (Fig.~\ref{fig:pitcheigs}) quantifies these insights and clearly demonstrate the effect of coarse-graining on the system (i.e., increasing $t_{max}$ while keeping the number of samples constant).
At very short time scales ($t_{max}<.05$), $x_0$ and the highest order $\alpha$ term are the main participants of the leading eigenvector, and $x_0$ soon falls off as $t_{max}$ increases; recall from Fig.~\ref{fig:pitchtraj} that this high-order term was equivalently able to bend the trajectory significantly until $t \approx 1$.
Around $t_{max}=10$, $r$ begins to have a noticeable influence on the observed trajectory, and correspondingly this is the point where $r$ becomes the dominant participant in the leading eigenvector.
For large $t_{max}$, the leading eigenvalue increases while all other eigenvalues decrease, indicating that the system's bifurcation is codimension one.
Note that in this range, small changes to the initial value $x_0$ have fallen all the way to the last eigenvector, indicating that the system has been allowed to run long enough that transient dynamics are removed, or at least have orders of magnitude less influence than any of the nuisance parameter $\alpha_i$'s. 
There is no significance to the fact that in this and subsequent ``rainbow diagrams'', the leading eigenvalue eventually begins to increase; this is simple case of an increasing line overtaking non-increasing ones and nothing inherent about the highest eigenvalue at small time scales. This can be confirmed by the change in color, indicating that the parameter responsible for the leading eigenvector has changed.

Similar figures can be produced for the saddle-node, transcritical, and subcritical pitchfork bifurcation classes. In each case, the eigenanalysis of the FIM indicates 
\begin{itemize}
    \item how long the system should be simulated, by the time it takes for the effect of the initial conditions to reach the least relevant eigenvector
    \item the codimension of the bifurcation, by the number of non-decreasing eigenvalues ($=1$ for each normal form),
    \item the participation factor of each parameter in the hyper/relevant directions by the square of the corresponding eigenvectors (asymptotically approaching 100\% $r$ in each normal form)
    \item the null-space of the bifurcation surface, making it possible to track the bifurcation hypersurface through parameter space.
\end{itemize}

These are relatively simple bifurcations, where the separatrix is the hyper-plane $r=0$. In more complicated situations where the separatrix is a nonlinear combination of bare parameters, this analysis identifies the vector normal to the separatrix. In principle, this local characterization could be extended to map that separatrix (along the hyper/relevant directions) through the high-dimensional parameter space.

%\begin{figure}[b]
%\includegraphics[scale=.4]{hopf_traj.png}
%\caption{\label{fig:hopftraj} The Hopf bifurcation, starting with $\mu=0.1$, $\beta_1=-1$ and $\beta_2=1$, and $y_0=(.03+\sqrt{0.1},0)$ to start just outside the limit cycle (marked with a *). The trajectory is similar for small changes to $\mu$ and the $\beta_i$'s, but the resulting endpoints (marked with circles) move much further for $\beta$ than for $\mu$.}
%\end{figure}

\begin{figure}[b]
    \centering
    \includegraphics[width=.98\linewidth]{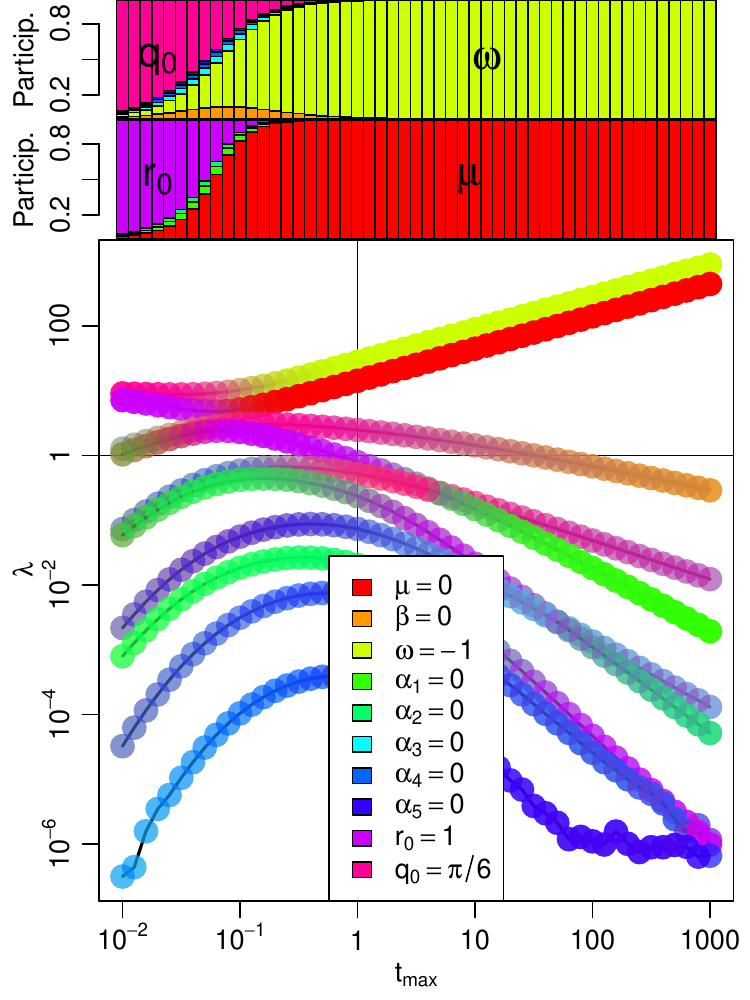}
    \caption{TWIG analysis of the Hopf bifurcation. The first of the hyperrelevant (rising) eigenvalues comes from the periodicity of the trajectory, whose velocity is set by $\omega$. The second hyperrelevant eigenvalue comes from the bifurcation itself, indicating that the Hopf bifurcation is codimension-1, and the bifurcation depends simply on $\mu$, and not some complicated combination of parameters. Note that the Hopf bifurcation is far easier to simulate at long time scales in polar form than in Cartesian coordinates.}
%    \caption{THIS IS FOR THE CARTESIAN HOPF: The eigenvalues of the Hopf bifurcation's FIM show a dependence on $\mu$ that switches from the leading to the second eigenvector after one rotation, and a strong but equal participation both $\beta$ parameters in the leading and third eigenvector (the $\beta$ parameters have opposite signs in the first vector, and the same sign in the third.)}
    \label{fig:hopfeig}
\end{figure}

Hopf bifurcations present more of a challenge: the limit cycle that emerges from a fixed point can have both its radius \emph{and} velocity manipulated by model parameters, which can potentially provide the false impression that the system has codimension 2, rather than the actual codimesion of 1. 
%Consider the following Hopf bifurcation in Cartesian coordinates, where, as above, additional high order terms have been added: \begin{equation}
%    \begin{aligned}
%        \dot x&=\beta_1y + \rho x + \alpha_1 \rho^2 x + \alpha_3 \rho y \\
%        \dot y&=\beta_2 y + \rho y + \alpha_2 \rho^2 y + \alpha_4 \rho x
%    \end{aligned}
%    \label{eq:hopf}
%\end{equation} where $\rho = \mu - (x^2+y^2)$. In this system, a limit cycle emerges of radius $\sqrt{\mu}$ when $\mu>0$. As shown in Fig.~\ref{fig:hopftraj}, trajectories beginning near the limit cycle quickly spiral onto it in a counter-clockwise fashion. Adjusting $\mu$ slightly (blue) increases the radius of the limit cycle, but the endpoint is at the same angle as the base case and therefore not very distant from the black endpoint. By contrast, increasing the $\beta$'s (red) increases the velocity of the cycle, so the red endpoint may well be far from the black endpoint if the dynamic system runs for many cycles. As before, the $\alpha_i$ parameters have only transient effects, and are not shown in Fig.~\ref{fig:hopftraj} for clarity.
%Actually, let's not worry about cartesian form.
Consider the following Hopf bifurcation in polar coordinates, where, as above, additional high order terms have been added: \begin{equation}
    \begin{aligned}
      \dot{r} &= \mu r - r^3 + \alpha_1 r^4 + \alpha_2 r^5  \\
      \dot\theta &= \omega + \beta r^2 +\alpha_3 r^3+\alpha_4 r^4 + \alpha_5 r^5
    \end{aligned}
\end{equation}
At the bifurcation point $\mu=0$, a fixed point at the origin expands into a limit cycle. The velocity of trajectories around this cycle are primarily driven by $\omega$, provided $y$ values are small. Note that the periodicity of the Hopf bifurcation introduces a second hyperrelevance to long-term dynamics. Infinitesimal changes to velocity make little difference to the final position of the trajectory $F(t_{max}; y,\theta)$ \emph{if $t_{max}$ is small}, but will have an increasing effect as $t_{max}$ grows. By contrast, $\mu$ is hyperrelevant because it is the bifurcation parameter. The increasing importance of these two parameters, in contrast to all others, is clearly illustrated in Fig.~\ref{fig:hopfeig}.

%This can also be quantified using the FIM eigenanalysis applied to the other bifurcation classes above. In Fig.~\ref{fig:hopfeig}, we see up until $t_{max} \approx 2\pi$, $\mu$ is only slightly more relevant than the $\beta$ parameters, but after one complete revolution, the trajectory has come very close to the limit cycle whose radius is set by $\mu$, and so changes to $\mu$ will not become any more relevant than they already are.
%\mkt{I do not think this statement is true.}
%By contrast, the velocity of the cycle, set by the $\beta$ values, continues to increase in relevance as every sample point in the cycle moves in response to these numbers. The relative importance of $\beta_1$ and $\beta_2$ is exactly split on both the first and third (not pictured) eigenvectors, while $\mu$ takes sole possession of the second eigenvector. Thus, the first eigenvector is hyperrelevant because it controls angular velocity, and the third is relevant because it controls absolute phase. Only the second is relevant due to its involvement in the bifurcation. 

\begin{figure*}
  \includegraphics[width=0.49\linewidth]{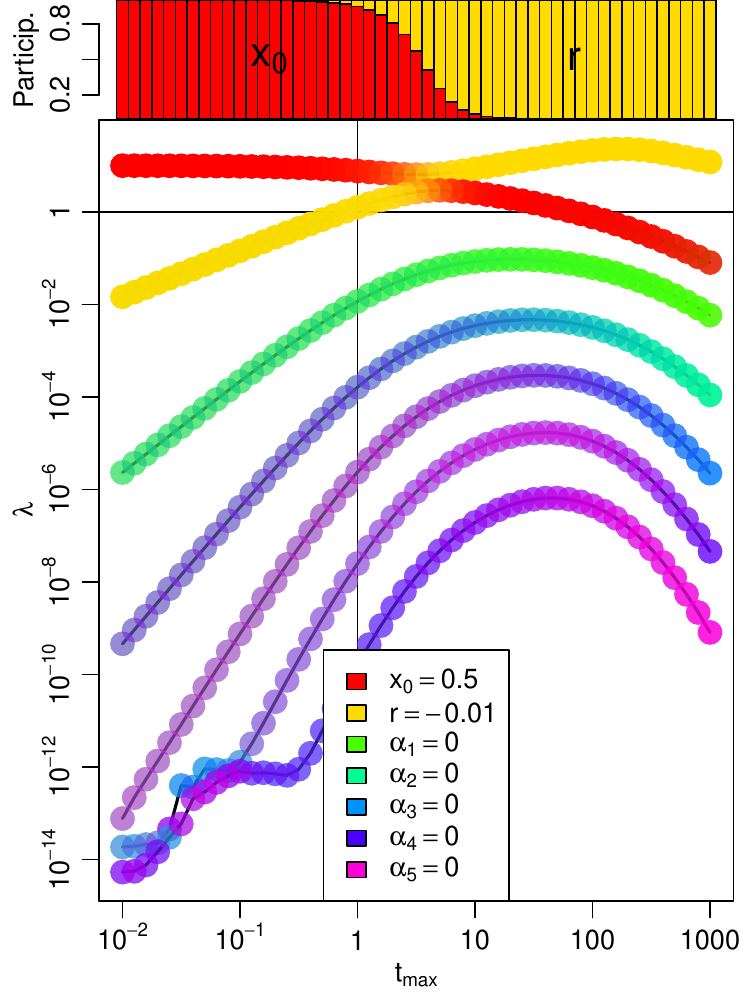}%
  \hfill
  \includegraphics[width=0.49\linewidth]{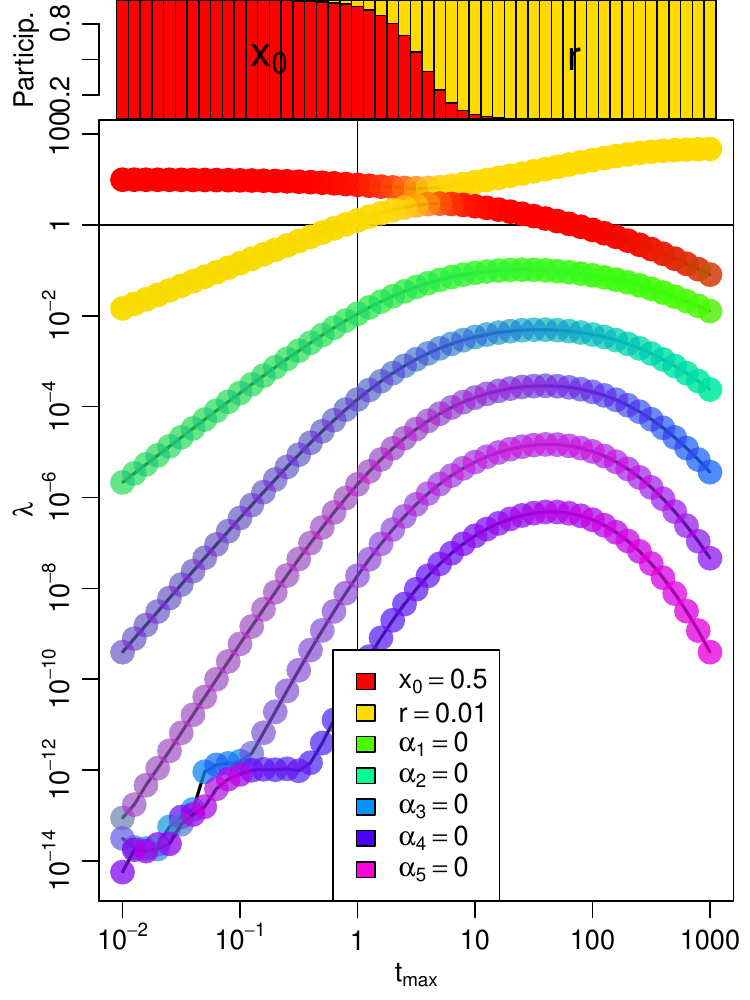}%
    \caption{TWIG analysis near-but-not-at the bifurcation values show the diagnostic pattern of an increasing eigenvalue at intermediate time scales, but then a decrease (in A) or an asymptote (in B; top yellow curve in the $t_{max}$ =100 to 1000 range in both). It is still relatively easy to identify bifurcation parameters and the bifurcation's codimension, though the fall-off becomes more pronounced, and thus the clarity of the analysis more obscured, the further one moves away from the bifurcation in either direction.}
    \label{fig:jonim}
\end{figure*}

As noted above, this ability to characterize all normal-form bifurcations depends on the ability to isolate changes in information due to the bifurcation itself.
% If the information does not change--say, because trajectories on both sides of an unstable point diverge--TWIG will not detect the bifurcation.
% In all likelihood, it will also not be possible to run simulations long enough to create figures like those in the paper before computer overflow occurs.
This depends on the only source of variation in long-term behavior coming from the bifurcation, and so the preceding analyses were conducted for systems exactly at the bifurcation point.
We now consider how the picture changes for dynamics near, but not exactly at, the bifurcation point.
Applying TWIG just to the left and right of the bifurcation point of a pitchfork ($r=\pm.01$) shows characteristic patterns (Fig.~\ref{fig:jonim}).
In these cases, we find that the bifurcation parameter is hyper-relevant on intermediate time scales ($10 < t_{max} < 100$ in Fig.~\ref{fig:jonim}).
However, on longer time scales ($t_{max} > 100$), the leading eigenvalue either decreases (Fig.~\ref{fig:jonim}A) or asymptotes (Fig.~\ref{fig:jonim}B) once the trajectories have converged to the fixed point, depending on whether the location of the fixed point can or cannot be shifted by changing the parameter values, respectively.
In other words, when approached from the $r<0$ side, small changes to $r$ don't move the equilibrium ($y(t)\rightarrow 0$ as $t\rightarrow\infty$), meaning the exact value of $r$ is irrelevant. But approaching from the $r>0$ side causes trajectories to run to $y(t)\rightarrow \pm \sqrt{r}$, meaning $r$ is relevant. 
Moving the system closer to bifurcation, this intermediate regime extends further and further, until at $r=0$ it occupies the entire trajectory and $r$ is hyperrelevant at all times $t>10$.

\begin{figure}
    \includegraphics[width=.98\linewidth]{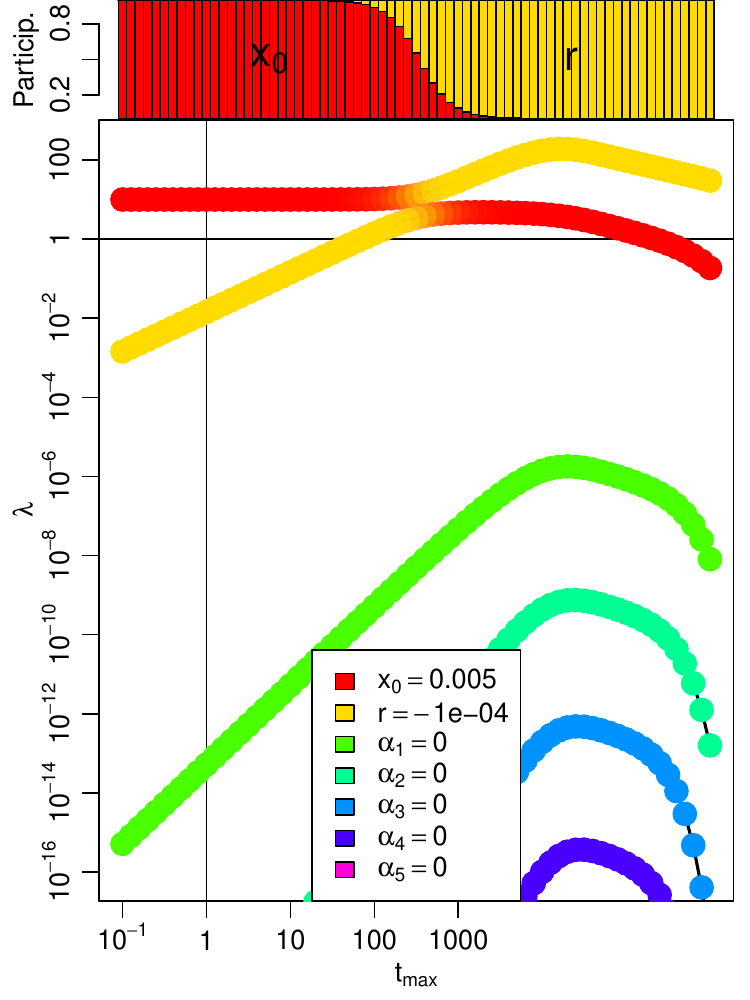}
    \caption{The subcritical pitchfork cannot be analyzed using TWIG at the bifurcation point ($r=0$) because the system is unstable. However, simulations slightly to the stable side of the bifurcation ($r\rightarrow 0^-$) reveals the bifurcation parameter, though because analysis happens off the bifurcation, the peak occurs at intermediate values instead of reaching an asymptote.}
    \label{fig:pfsub}
\end{figure}

In general, being slightly off the bifurcation obscures the effect of the bifurcation parameter to an extent proportional to the distance from the bifurcation.
This is particularly useful in the case of hemi-stable bifurcations, which need to be approached from the stable side or else test trajectories will diverge to infinity (and cause computer overflow). 
In the case of the subcritical pitchfork, at the bifurcation itself ($r=0$) the system is unstable. 
However, at values of $r\rightarrow0^-$, just less than bifurcation value, TWIG can be performed and the bifurcation characterized as above (Fig.~\ref{fig:pfsub}).

% On the fixed-point (negative) side, there are intermediate time scales on which $r$'s importance is clear, but on long time scales 
% Similarly, on the limit-cycle (positive) side, the importance of $r$ begins to tail off, albeit on longer time scales, for similar reasons.
% In both cases, $r$ is clearly the driver of the leading eigenvalue at long scales.

\section{\label{sec:abnorm}Bifurcations in Non-normal Forms}
Equations describing real systems are not typically written in one of these normal forms. So even when a researcher knows a system contains a bifurcation, it might not be apparent which one of these it is. For example, a model of a bead on a rotating hoop $$mr \frac{\partial^2\phi}{\partial t^2}=-b~\frac{\partial \phi}{\partial t}-mg \sin \phi + mr \omega^2 \sin \phi \cos \phi$$ has a supercritical pitchfork bifurcation, though it might require simulating many values of $r$ and $\omega$ to appreciate this~\cite{Strogatz2015}. 
%[The system is described in section 3.5, but has the same form as an overdamped pendulum in 4.4 from ]; won't appear in refs because Strogatz already has a comment from line 384?
Similarly, the equation: 
\begin{equation}
  \label{eq:nontrivialtranscritical}
  \dot{x}=r\ln x +x-1+\alpha_1 (x-1)^2+ \alpha_2 (x-1)^3+\cdots
\end{equation}
contains a transcritical bifurcation at $x=1$ when $r=-1$. However, this only becomes clear after reparameterizing the equation by $R=r+1$, and $X=\frac{r}{2}(x-1)$, when the equation assumes the normal form $\dot{X}=R X-X^2+\mathcal{O}(X^3)$. Such a substitution might not be immediately apparent to a researcher; however, time-widening information geometry clarifies the situation. 

\begin{figure}
    \centering
    \includegraphics[width=.98\linewidth]{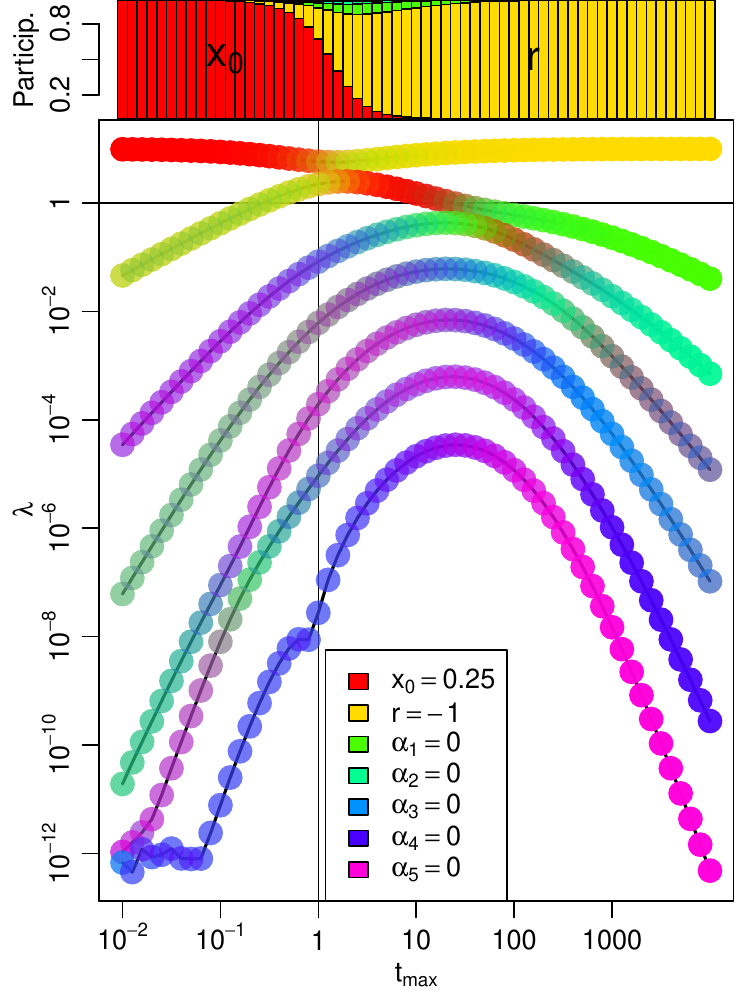}
    \caption{Equations such as Eq.~\eqref{eq:nontrivialtranscritical} that are not in normal form can be interpreted using the same procedure as for normal form bifurcations. As above, the presence of just one non-decreasing eigenvalue, whose corresponding eigenvector is dominated by the single parameter $r$, indicates that the system has codimension 1 and the bifurcation parameter involves only $r$. The relevant (not hyperrelevant) leading eigenvalue is characteristic of a transcritical bifurcation.}
    \label{fig:easy2}
\end{figure}

If the dynamics in Eq.~\eqref{eq:nontrivialtranscritical} are run long enough, we observe that one eigenvalue is relevant while all others are irrelevant.
Furthermore, the corresponding participation factor becomes dominated exclusively by $r$ (Fig.~\ref{fig:easy2}).
This tells us that (1) the process has codimension 1, and (2) the reparameterization involves only $r$.
We confirm that our analysis has converged since the initial condition $y_0$ is the dominant participation factor in the smallest eigenvalues.
However, we note that convergence occurs at a somewhat larger value of $t_{max}$ than in the normal form examples above. Also note that transcritical bifurcations have a leading eigenvalue that is \emph{relevant} rather than hyperrelevant, due to a quirk of the normal-form algebra. See Appendix \ref{ap:trans} for a thorough explanation.

\begin{figure}
    \centering
    \includegraphics[width=.98\linewidth]{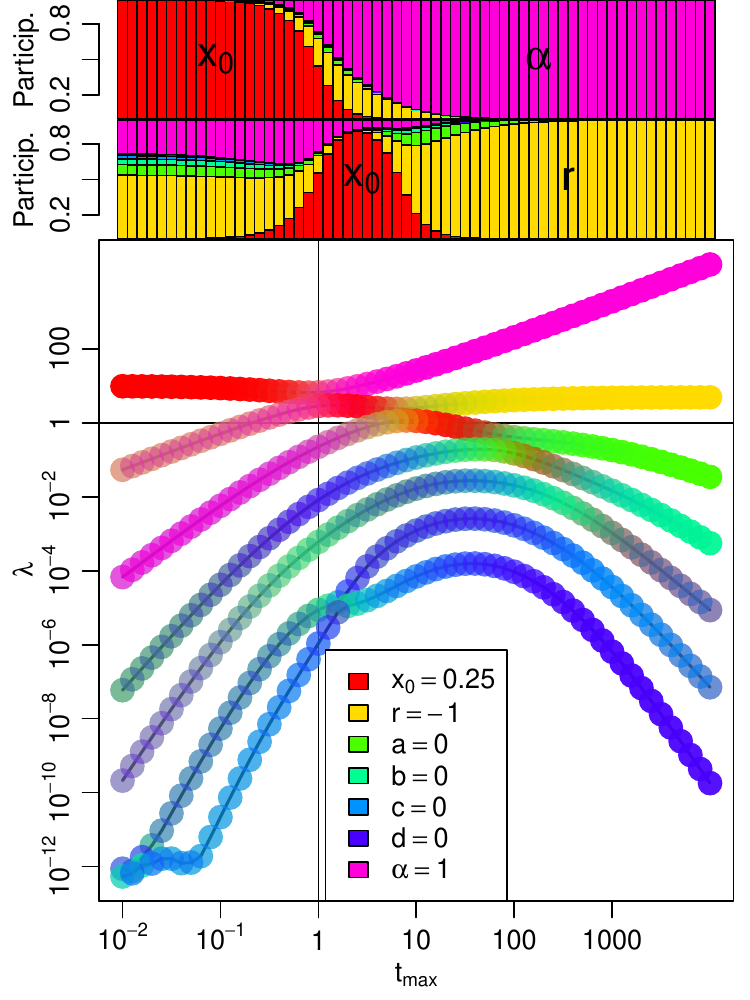}
    \caption{A difficult non-normal-form transcritical bifurcation such as Eq.~\eqref{eq:modifiedtranscritical} can be extremely challenging to analyze analytically, but sloppy analysis indicates one hyperrelevant parameter (corresponding in this case to a saddle-node) and one relevant parameter (as usual, indicating transcritical bifurcation). This means that this system has a bifurcation of codimension two. Note that the participation factor of the two leading eigendirections runs to 1.0 in the direction of $\alpha$ and $r$ respectively, indicating that the system can be placed into normal form without a complicated recombination of parameters.}
    \label{fig:hardbif}
\end{figure}

But what happens when the situation is not so straightforward? Modifying the above example to the equation
\begin{equation}
  \label{eq:modifiedtranscritical}
  \dot x =r \ln(x)+ a (x-\alpha) +b(x-\alpha)^2+\cdots
\end{equation}
should still have a transcritical bifurcation for certain parameter values, but no simple reparameterization to create a normal form exists. From above, we can recognize that when a transcritical bifurcation occurs at $x=1$ for $r=-1, \alpha=1$. However, when $\alpha \neq 1$, in the neighborhood of $x=\alpha$ all the power terms are zero, but the term $r \ln(x) > 0$ if $\alpha<1$, suggesting that no fixed point exists in that region. 
The appearance or disappearance of a fixed point is the hallmark of a saddle-node bifurcation, and indicates that allowing a bit of variability in the fixed point's location has introduced a second codimension to the dynamical system.
This is borne out by TWIG analysis (Fig.~\ref{fig:hardbif}), which shows that the equation indeed produces a hyperrelevant eigenvector corresponding to the saddle-node parameter $\alpha$, which controls the existence--not just the location--of an equilibrium. The transcritical bifurcation still exists and is controlled by $r$, as implied by the previous analysis. This example shows that even in situations with two different bifurcation classes, neither of which can be reparameterized into normal form, TWIG still allows us to efficiently and unambiguously identify co-dimension and bifurcation parameters.

%The concept of codimension may require some clarification here. It refers to the number of parameters (or rather, the minimum number when optimally reparameterized) necessary to tune in order to achieve a change in behavior, or nonhomeomorphic topological difference.  Consider eqn \ref{eq:hopf}, illustrating the Hopf bifurcation. Recall when $\mu \leq 0$, any trajectory spirals into the origin, and when $\mu>0$ to a limit cycle of radius $\sqrt{\mu}$ around the origin. Here, $\beta_1$ and $\beta_2$ control the velocity and direction of the limit cycle (set to {1,-1} in normal form), and the $\alpha$'s are nuisance parameters (0 in normal form). The system has two dimensions in phase space, which could also be data space if $(x,y)$ are both observed, and/or prediction space if forecasts are made of both $x$ and $y$. It has seven dimensions in parameter space--not eight, since $r$ is determined by the system state. However, the non-homeomorphic change in topology from a fixed point to a limit-cycle is controlled solely by $\mu$, and therefore the bifurcation's codimension is 1, meaning the dimension of the "bifurcation surface" (the hypersurface separating the basin of fixed point behavior from the basin of limit cycle behavior) is six. Not conflating these four different kinds of "dimension" is one of the significant challenges of bifurcation analysis. 

\subsection{A biophysical example}

Glycolysis is a multi-step process which uses the bond energy of glucose to catabolize energy-carrying biomolecules easily usable by cells, which represents one of the dominant processes of all heterotrophic life on earth. A bottleneck in this crucial process is the phosphorylation of fructose-6-phosphate into fructose-1,6-bisphosphate catalyzed by the enzyme phosphofructokinase. The complicated five-species mass-action equation describing this reaction's kinetics can be simplified using Tikhonov's theorem and assuming low concentrations of ATP to the simple dimensionless system:~\cite{Selkov1968,Tikhonov1948}
\begin{equation}
\begin{aligned}
\dot{x}&=-x+ay + c_1 x^2 y +c_2 x^3\\
\dot{y}&=b - ay + c_3x^2 y + c_4 y^2
\end{aligned}
\label{eq:glyco}
\end{equation}
where $x$ and $y$ are the concentrations of ADP and F6P respectively, and the four $c_i$ constants are nuisance parameters added to mask the system dynamics. There is a curved bifurcation surface that separates the range of kinetic parameters $a,b$ which lead to either a fixed point at $(b,\frac{b}{a+b^2})$ when $c_1=1,~c_3=-1$ as in the canonical model, or a stable limit cycle. The separation between the fixed point and limit cycle regimes has the form $b^2=\frac 1 2 \left(1-2a\pm\sqrt{1-8a}\right)$~\cite{[Figure 7.3.7 in ]Strogatz2015}. 
The resulting oscillations in glycolytic activity predicted by this analysis have been observed \emph{in vivo} since the early 1970s~\cite{Chance1973}.

\begin{figure}
    \centering
    \includegraphics[width=.98\linewidth]{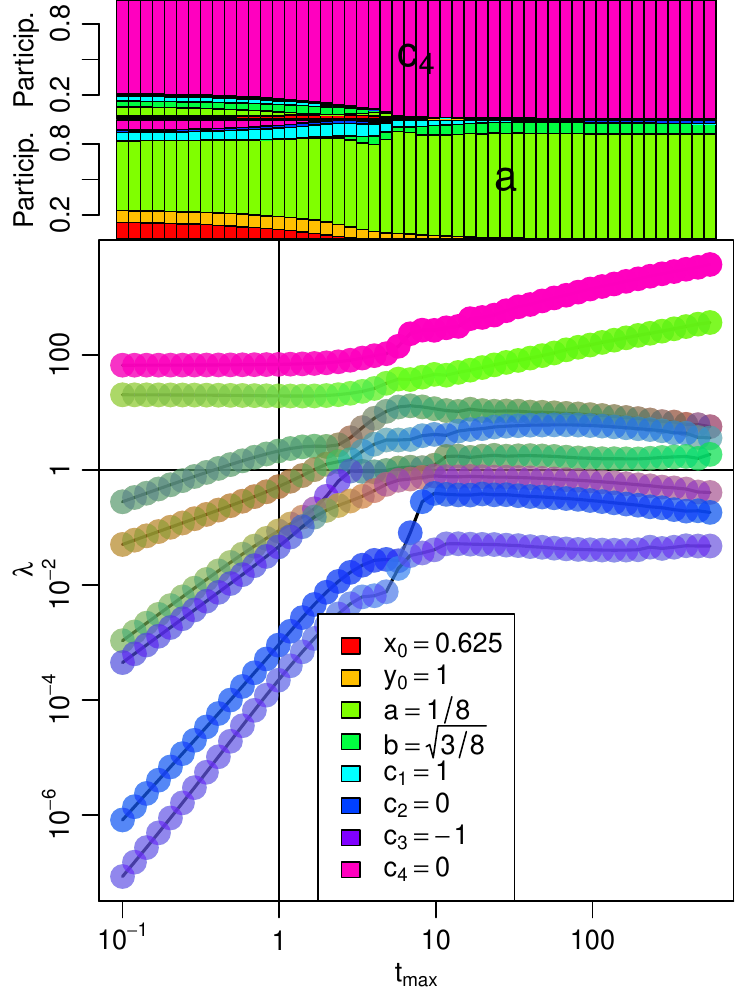}
    \caption{Analysis of the ``glycoscillator" bifurcation (Eq.~\ref{eq:glyco}). The frequency of the oscillations are driven by $c_4$, while the radius of oscillations can be controlled with just one of the $a,b$ parameters discovered by Sel'kov~\cite{Selkov1968}. }
    \label{fig:glycofig}
\end{figure}

A TWIG analysis of this system provides several insights, summarized in Fig.~\ref{fig:glycofig}. First, even though the separatrix between fixed point and limit cycle in $a,b-$space is a nonlinear curve, because $b$ can be reparameterized as a function of $a$, it is codimension one. Second, the ``nuisance" parameter $c_4$ introduces a change in the period of the oscillations, which means that infinitesimal changes in its value cause larger deviations in final trajectory the longer the simulation runs. This shows up as a hyperrelevant direction in TWIG; however, as discussed above, it is \emph{not} a second codimension.

\section{\label{sec:chaos}Chaotic systems}

Systems showing chaotic behavior have long represented a challenge to traditional categories of thinking, and the difficulty of distinguishing deterministic chaos from randomness is practically its own subdiscipline~\cite{Briggs1989, Holmes1990, Strogatz2015, Sugihara1990a, Hsieh2008}. In the context of TWIG analysis, there are two characteristics of the system that need to be considered carefully.

First, unlike other systems considered here, one hallmark of chaos is long-term sensitive dependence on initial conditions, or the ``butterfly effect''. Because of this, a TWIG analysis carried out in the chaotic regime, in contrast to Fig.~\ref{fig:pitcheigs} where the parameter $x_0$ becomes the least relevant, will classify initial condition parameters as relevant. Note that if the chaotic system produces a strange attractor, then the initial conditions will change the location of the system on the attractor at long time scales, but not the shape of the attractor itself, which prevents these parameters from becoming hyperrelevant. That is, the maximum distance between two trajectories begun at slightly different initial conditions will eventually saturate on opposite sides of the attractor, and not increase without bound.

Second, the four classic examples of chaotic systems approach chaos through a complicated series of bifurcations, rather than a singular event as in the normal-form bifurcations above. The logistic map famously contains a period-doubling ``bifurcation cascade'', with the distance between these bifurcation events decreasing geometrically by the Feigenbaum constant $\alpha$ universally~\cite{Feigenbaum1978, Feigenbaum1979, [see also Chapter 10.6 of ]Strogatz2015}. There is thus a ``fuzzy boundary'' between the periodic behavior of, say, an 8-cycle and the chaotic region as we pass through the increasingly narrow 16-cycle region, 32-cycle region, and so forth. The H\'{e}non map experiences a similar bifurcation cascade along the line $b=0.3$ as $a$ increases from 1 to 1.5~\cite{Hitzl1985}, while the R\"{o}ssler attractor has a bifurcation cascade in the opposite direction on the plane $a=0.2,~c=5.7$ as $b$ decreases from 1.5 towards 0~\cite{Rossler1976}. As we show below for a R\"{o}ssler system, these boundaries are not just fuzzy, but also fractal. Most complex of all, the Lorenz system experiences a pitchfork bifurcation at $r=1$, whose two stable points then experience Hopf bifurcations at $r \approx 24.74$, while the unstable point undergoes a ``homoclinic explosion'' at $r \approx 13.926$ that produces an ``a thicket of infinitely many saddle-cycles and aperiodic orbits \cite{Strogatz2015}.'' If even these pedagogical ``toy models'' of chaos have such indeterminate boundaries, it is likely that examples of chaotic systems encountered ``in the wild'' will as well.

While the FIM may be evaluated at any point in this fuzzy region between order and chaos, its interpretation is less clear.
The eigenvectors, which indicate the direction normal to the separatrix in other systems, lose this meaning since there is no direction normal to a fractal surface.
Note that this also holds true for the intermittency route to chaos as well. Abrupt changes to chaos, with or without smooth changes in fractal dimension, also exist and would be expected to give cleaner results in the TWIG analysis~\cite{Bleher1989, Lai1999}, but unfortunately are expected to be less common and less familiar to readers.

\begin{figure}
    \centering
    \includegraphics[width=0.98\linewidth]{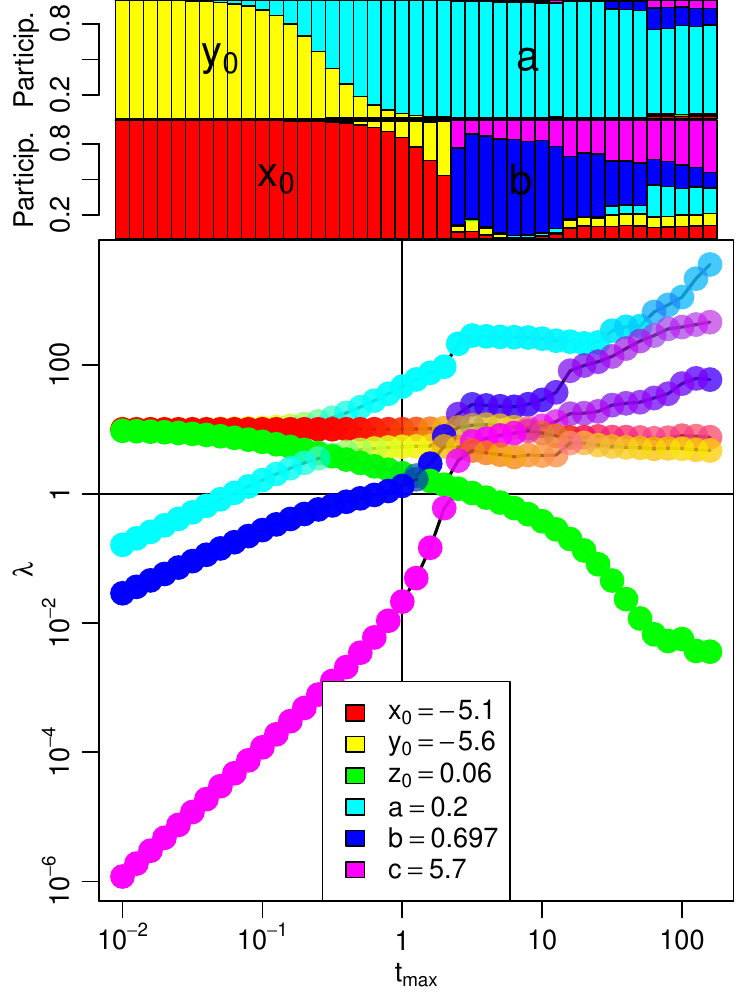}
    \caption{TWIG analysis of the R\"{o}ssler attractor, a chaotic system, evaluated in the region of rapid period doubling just before the onset of chaos. Due to the butterfly effect, the initial conditions remain relevant at long time scales, and cannot be used to determine appropriate simulation length. However, excluding these from analysis, we are still able to qualitatively see that there is one hyperrelevant direction, dominated by $a$. This came as a surprise to the authors, because the bifurcation region was approached by changing values of $b$ until a period-doubling cascade was observed, yet TWIG uncovered a greater sensitivity to $a$ than $b$ even in this region. This was confirmed by sampling the parameter space in Fig.~\ref{fig:ross_parsp}.}
    \label{fig:ross_eig}
\end{figure}

\begin{figure}
    \centering
    \includegraphics[width=0.98\linewidth]{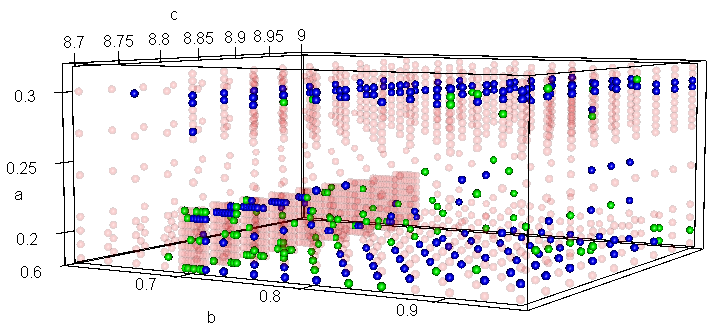}
    \caption{The parameter space in the period-doubling region of the R\"{o}ssler attractor shows flat sheets of 8-cycle behavior (solid blue spheres) sandwiched between chaos (transparent red) in the $a$ direction. Green spheres are simulations difficult to classify as either 8-cycle or chaotic.}
    \label{fig:ross_parsp}
\end{figure}

That being said, TWIG still can provide powerful qualitative insights into the nature of a period-doubling chaotic system. The R\"{o}ssler attractor defined by 
\begin{equation}
    \label{eq:rossler}
    \begin{aligned}
        \dot{x} &= -y -z\\
        \dot{y} &= x + a y\\
        \dot{z} &= b + z(x-c)
    \end{aligned}
\end{equation}
has a well-known period doubling map revealed by decreasing $b$ along the parameter-space line $a=.2, c=5.7$, with the fuzzy transition from an 8-cycle to chaos occurring in the region near $b\approx .70$. TWIG analysis carried out near this point reveals that while changes to $b$ or $c$ in this region can lead to long term divergent behavior, changes to $a$ have a much stronger effect (Fig.~\ref{fig:ross_eig}). In other words, even though we ``walked up'' to the bifurcation region in the $b$ direction, TWIG was able to tell us that the fuzzy bifurcation boundary was strongly angled normal to the $a$ direction. This insight is not found in the usual treatments of the R\"{o}ssler attractor~\cite{Rossler1976, Guckenheimer1983, Strogatz2015}, but can be easily verified by simulating the system in Eq.~\ref{eq:rossler} at many sample parameter values in the region around the bifurcation. This reveals flat ``sheets'' of periodic behavior sandwiched between strata of chaos in the $a$ direction (Fig.~\ref{fig:ross_parsp}); these sheets can eventually be encountered for a fixed value of $a$ by moving far enough in $b$ or $c$, which is essentially the process diagrammed in the period-doubling map with which we started this exercise. 

Above, we made the claim that TWIG analysis could be used to determine four characteristics of the system, the first being the length of time to run an analysis by the decay of sensitivity to initial conditions. For chaotic systems, this is no longer the case due to the butterfly effect. However, by removing the initial values as parameters, we see that TWIG can still be used \emph{qualitatively} to determine the other three characteristics: bifurcation co-dimension, the null space, and the (fuzzy) normal to the bifurcation region. 

\section{\label{sec:conc}Conclusion}

Gradual time-dilation of the Fisher information matrix as realized by our Time-Widening Information Geometry (TWIG) analysis is an efficient way of characterizing bifurcations in a dynamical system. Researchers have long used eigenanalysis of $\mathcal I$ to characterize the ``sloppiness" of a system, \emph{i.e.} its exponential range of sensitivities to parameter changes, and recently leveraged this accumulated expertise with coarse-graining to understand phenomena occurring at distinct time scales~\cite{Raju2018,Chachra2012}. Building on these insights, we here demonstrate that as $t_{max}$ increases, the changing eigenvalues of $\mathcal I$ (and the composition of the corrresponding eigenvectors) allow us to (1) characterize the codimension of the bifurcation, (2) quantify the participation of each bare parameter in the bifurcation, (3) characterize the bifurcation's hyper-surface, and (4) have an internal check on the length of time necessary to simulate the system to reach equilibrium. These are substantial insights to be gained relatively cheaply. Sloppy bifurcation analysis constitutes a powerful tool to supplement traditional analytical analysis \cite{Rasband1990,Strogatz2015}, and other specialized analytical tools for high-dimensional problems~\cite{Rand,Feigenbaum1978,Bettonvil1997,Waldherr2007,Guckenheimer1983,Kirk2008,Gul2015}. 

Insights derived from TWIG are useful not just for theoreticians interested in characterizing a bifurcation or reparameterizing a system to emphasize the bifurcation; it is also critical for the process of fitting parameter values. The rainbow plots in this paper demonstrate that at the bifurcation point, simulations frequently show a separation of over 10 orders of magnitude in their parameter sensitivity, a gap that gets larger the longer the simulations run (or the more data is collected, in an experimental context). 
If researchers care about fitting all parameters, it is crucial to recognize that the effect of hyperrelevant parameters will overwhelm the others, so only if these parameters are fixed in the experiment can the less relevant ones be inferred~\cite{Kirk2008, Roesch2019}.
Future work may naturally extend the method to large systems including those derived from partial differential equations.

Our TWIG analysis has some inherent limitations.
It presupposes that the model can be simulated on at least one side of the bifurcation to arbitrarily long times, \emph{i.e.} it analyzes stable dynamics on the threshold of instability.
A bifurcation that switches between two different forms of instability will not be easily detectable with this method, as trajectories will diverge on both sides of the bifurcation.
However, such doubly-unstable bifurcations may be of limited practical interest anyway, as loss of stability is generally a far more common real-world problem than a change in the type of instability of a system that never was stable to begin with.
Hemi-stable points (as in saddle-node or subcritical pitchfork bifurcations) are easily analyzed when approached from the stable side (see Fig.~\ref{fig:pfsub}); otherwise test trajectories can diverge beyond computer tolerance at moderate time scales.
A notable limitation of the method as presented here is the inability to analyze hyperbolic fixed points.
Future work may additionally leverage center manifold techniques to investigate bifurcations in such systems.
We note here that absolutely unstable fixed points (\emph{i.e.}, where every eigenvalue is positive) can be conveniently analyzed in TWIG simply by running time backwards, and analyzing trajectories at ever-closer instants to the initial divergence from the instability.

Because it is a particularly efficient method of determining important information about high-dimensional bifurcations, we anticipate that TWIG will be useful in situations with many components where one or a few bifurcations are expected in each component. These include power grids, circuit boards, interatomic models, complex protein regulatory networks, and ecosystem-based management systems of multiple interacting populations. Such complexity presents substantial difficulties for closed-form analysis but can be tamed with insights gleaned from this method.

\emph{Acknowledgements:} We thank Archishman Raju for helpful discussions and comments on the manuscript. 
We also thank two anonymous reviewers for their careful reading and suggestions for extensions of TWIG.
This work was supported by the US National Science Foundation under Award NSF-1753357.

\appendix

\section{\label{apx:SN}FIM of Saddle-Node Bifurcations}

The normal form of the saddle-node bifurcation is 
\begin{equation}
    \frac{dy}{dt} = r+y(t)^2 +\alpha_1 y(t)^3 +\alpha_2 y(t)^4 + \dots
\end{equation}
This differential equation can be solved locally when all parameters $\overrightarrow{\theta}=0$, which happens to be the bifurcation point of the system. At that point: 
$$\frac{dy}{dt} = y^2 \rightarrow\frac{dy}{y^2} = dt$$ 
Integrating both sides yields
\begin{align*}
-\frac{1}{y}\Big|_{y_0}^{y(t)} &=t\Big|_{0}^{t} \\ 
\frac{1}{y_0}-\frac{1}{y(t)}&=t \\ 
y(t)=\frac{y_0}{1-y_0 t} \numberthis
\end{align*}
This implies there is a singularity at $t=1/y_0$, so a proper coarse-graining procedure will involve taking data from $t$ = 0 to some value near $1/y_0$, say $0.99/y_0$. 
We avoided this singularity by using negative values for $y_0$ and were therefore able to run simulations to large values of $t_{max}$.
As noted in Eq.~\ref{eq:JTJ}, to find the FIM of a system it is only necessary to find the Jacobian, so we need only find the first partial derivative of these data with respect to each parameter in the model.

\subsection{Partial derivative of $r$}
Let the $\alpha_i$'s=0. The derivative of the normal form w.r.t. $r$ becomes:

\begin{eqnarray}
\frac{\partial}{\partial r}\left(\frac{\partial y}{\partial t} = r+y^2 \right)\nonumber \\
\frac{\partial ^2 y}{\partial r \partial t}=1+2y\frac{\partial y}{\partial r}
\end{eqnarray}

We let $w=\frac{\partial y}{\partial r}$, and this becomes $\frac{\partial w}{\partial t}=1+2yw$, which requires the use of an integration factor to solve~\cite{[See equation 14.2.1 in ]Hassani2013}. If $p_1 x'+p_0 x=q$ then
\begin{equation}
x=\frac{1}{\mu p_1}\left[C+\int\mu q dt \right] \text{where}~ \mu=p_1^{-1}\text{exp}\left(\int\frac{p_0}{p_1}dt \right)
\end{equation}

Allowing $p_1=1$, $p_0=-2y$, $q=1$ implies that
$$\begin{aligned}
\mu &= 1^{-1} \exp\left(\int\frac{-2y}{1}dt \right) \\
&= \exp\left(-\int \frac {2y_0 dt}{1-y_0t} \right) \\
&= \exp (2\ln(1-y_0t)) \\
&= (1-y_0t)^2
\end{aligned}
$$
Therefore,
$$\begin{aligned}
w &= \frac{C+\int(1-y_0t)^2 dt}{(1-y_0t)^2} \\
&=\frac{C-\frac{(1-y_0t)^3}{3y_0}|_0^t } {(1-y_0t)^2} \\
&=\frac{C+\frac{1-(1-y_0t)^3}{3y_0} } {(1-y_0t)^2}
\end{aligned} $$

Recall this function is being evaluated at the initial condition, where the partial derivative $w=\frac{\partial y}{\partial r} = 0$ (i.e., changes to $r$ do not change $y_0$). This implies that $C=-\frac{1-(1-y_0t)^3}{3y_0}$; when $t=0$ this further reduces to $C=0$. Therefore.
\begin{equation}
\frac{\partial y}{\partial r} = \frac{1-(1-y_0t)^3}{3y_0(1-y_0t)^2}
\end{equation}

\subsection{Partial derivative of $\alpha_1$}
Using the same procedure as above, 
\begin{align*}
\frac{\partial}{\partial \alpha_1}\left( \frac{\partial y}{\partial t} \right) &= \frac{\partial}{\partial \alpha_1}\left(y^2 + \alpha_1 y^3 \right) \\
\frac{\partial ^2 y}{\partial \alpha_1 \partial t} &= 2y\frac{\partial y}{\partial \alpha_1}+y^3+\cancel{3y^2\alpha_1\frac{\partial y}{\partial \alpha_1}} \\
\frac{\partial w}{\partial t}&=2yw+y^3 \numberthis
\end{align*}

Note on the second line, we are able to cancel the third term because we are evaluating the slope where $\alpha_1$ is zero. On the last line, note that $p_0$ and $p_1$ are the same as for $r$, so as above $\mu=(1-y_0t)^2$, but since now $q=y^3$:

\begin{align*}
w&=\frac{C+\int y^3(1-y_0t)^2 dt}{(1-y_0t)^2} \\
&=\frac{C+\int \left(\frac {y_0}{1-y_0t}\right)^3(1-y_0t)^2 dt}{(1-y_0t)^2}\\
&=\frac{C+\int \frac {y_0^3 dt}{1-y_0t} }{(1-y_0t)^2}\\
&=\frac{C-y_0^2\log{(1-y_0t)}|_0^t } {(1-y_0t)^2} \\
&=\frac{C-y_0^2\log{(1-y_0t)} } {(1-y_0t)^2} 
\end{align*}

Again, assuming $w=t=0 \rightarrow C=0$, so
\begin{equation}
\frac{\partial y}{\partial \alpha_1} = - \frac{y_0^2\log{(1-y_0t)} } {(1-y_0t)^2}    
\end{equation}

\subsection{Partial derivatives of higher-order $\alpha$'s}

Higher order terms in the series are of the form $\alpha_n y^{n+2}$ and so
\begin{equation}
\begin{aligned}
\frac{\partial}{\partial \alpha_n}\left(\frac{\partial y}{\partial t} = y^2 + \alpha_ny^{n+2} \right) \\
\frac{\partial ^2 y}{\partial \alpha_n \partial t}=2y\frac{\partial y}{\partial \alpha_n}+y^{n+2} + \cancel{(n+2)y^{n+1}\alpha_n\frac{\partial y}{\partial \alpha_n}} \\
\frac{\partial w}{\partial t}=2yw+y^{n+2}
\end{aligned}
\end{equation}

As above, we are able to cancel $(n+2)y^{n+1}\alpha_n\frac{\partial y}{\partial \alpha_n}$ because we are solving for slopes about the point $\alpha_n=0$. With the same value of $\mu$, we use integration factors to demonstrate:

$$
\begin{aligned}
w&=\frac{C+\int \left(\frac {y_0}{1-y_0t}\right)^{n+2}(1-y_0t)^2 dt}{(1-y_0t)^2}\\
&=\frac{C+\int \frac {y_0^{n+2} dt}{(1-y_0t)^n} }{(1-y_0t)^2}\\
&=\frac{C-\frac{y_0^{n+1}}{1-n} (1-y_0t)^{1-n} |_0^t } {(1-y_0t)^2}\\
&= \frac{C+\frac{y_0^{n+1}}{n-1} (1-(1-y_0t)^{1-n})  } {(1-y_0t)^2}
\end{aligned}
$$

Which again implies that $C=0$ at the initial condition $t=0$, and so for $n>1$ we can say
\begin{equation}
\frac{\partial y}{\partial \alpha_n} = \frac{y_0^{n+1} (1-(1-y_0t)^{1-n})  } {(1-n)(1-y_0t)^2}    
\end{equation}

Recall that the Jacobian of our system is
\begin{equation}
J=\begin{bmatrix}
\frac{\partial y_0}{\partial r} & \frac{\partial y_0}{\partial \alpha_1} &  \frac{\partial y_0}{\partial \alpha_2} & \dots \\
\frac{\partial y_1}{\partial r} & \frac{\partial y_1}{\partial \alpha_1} &  \frac{\partial y_1}{\partial \alpha_2} & \dots \\
\dots & \dots & \dots & \dots 
\end{bmatrix}    
\end{equation}

Because the Fisher information matrix $\mathcal{I}=J^T J$, we can see that element $\mathcal{I}_{1,1}=\left(\frac{\partial y}{\partial r}\right)^2$ will be $\mathcal{O}(t^2)$ because $\frac{\partial y}{\partial r}$ is $\mathcal{O}(t^1)$; all other elements will be a lower order of $t$. Thus, at long time scales, the FIM's element (1,1) will grow faster than all other elements, and therefore the most relevant parameter is clearly $r$. 

In the case where $\mathcal{I}$ is being derived from data (or from noise added to a non-/normal form equation), the importance of $r$ can be evaluated by increasing $\sigma^2 \propto y_0^{-3}$. Since, by the central limit theorem standard error $\sigma^2 \propto n$, then the number of time points sampled should decrease as $n \propto y_0^{-3}$.

\section{\label{ap:trans}FIM of Transcritical Bifurcations}

These have a similar normal form as the saddle-node bifuractions above:
$$\frac{dy}{dt} = ry(t)-y(t)^2 +\alpha_1 y(t)^3 +\alpha_2 y(t)^4 + ...$$
However, the change of sign in the second term causes the solution to the differential equation to also have a changed sign:
\begin{equation}
\begin{split}
\frac{dy}{dt} = -y^2 \rightarrow -\frac{dy}{y^2} = dt \rightarrow \frac{1}{y}\Big|_{y_0}^{y(t)}=t\Big|_{0}^{t} \\ \frac{1}{y(t)}-\frac{1}{y_0}=t \rightarrow y(t)=\frac{y_0}{1+y_0 t}
\end{split}
\end{equation}
Now the singularity occurs at $t=-\frac{1}{y_0}$, which generally only complicates the coarse-graining if initial conditions are negative. 

\subsection{Partial derivative of $r$}

The full solution to the partial derivative of $r$ is somewhat complicated because it depends on $y$:
\begin{equation}
\begin{split}
\frac{\partial}{\partial r}\left(\frac{\partial y}{\partial t} = ry-y^2 \right) \\
\frac{\partial ^2 y}{\partial r \partial t}=r\frac{\partial y}{\partial r}+y-2y\frac{\partial y}{\partial r} \\
\frac{\partial w}{\partial t}=w(r-2y)+y
\end{split}
\end{equation}
where $w=\frac{\partial y}{\partial r}$. Recall that the derivative is being evaluated where $r=0$, and so we can argue that

\begin{align*}
\frac{\partial w}{\partial t}+2yw & =  y \rightarrow \\
  \mu & =  \exp\left(\int\frac{2y_0dt}{1+ty_0}\right) \\
                                  & =  \exp[2\log(1+ty_0)] \\
  & =  (1+ty_0)^2 \numberthis
\end{align*}

Using our integration factors, we see:
\begin{align*}
    w&=\frac{C+\int{(1+ty_0)^2 \frac{y_0}{1+ty_0}dt} }{(1+ty_0)^2} \\
&=\frac{C+y_0t(1+\frac{y_0t}{2})} {(1+ty_0)^2} \rightarrow C=0\\
&=\frac{y_0t(2+y_0t)} {2(1+ty_0)^2} = \frac{\partial y}{\partial r} \numberthis
\end{align*}

Note that in the limit that $t\rightarrow\infty$, this expression is order 0 for $t$; therefore, unlike the other bifurcation classes, transcriticals are expected to have a \emph{relevant}, rather than a hyperrelevant, leading eigenvalue. This was confirmed with simulations (see Fig.~\ref{fig:easy2}).

\subsection{Partial derivative of $\alpha_1$}

The derivative can be set up as:
\begin{align*}
\frac{\partial}{\partial \alpha_1}\left(\frac{\partial y}{\partial t} = -y^2 + \alpha_1 y^3 \right) \\
\frac{\partial ^2 y}{\partial \alpha_1 \partial t}=-2y\frac{\partial y}{\partial \alpha_1}+3\alpha_1y^2\frac{\partial y}{\partial \alpha_1} +y^3 \\
\frac{\partial w}{\partial t}=-2yw+y^3 \numberthis
\end{align*}

Since we already know that $\mu=(1+ty_0)^2$, it follows that
\begin{align*}
w=& \frac{C+\int{(1+ty_0)^2 \left(\frac{y_0}{1+ty_0}\right)^3} }{(1+ty_0)^2} \\
=&\frac{C+y_0^2\int{\frac{y_0}{1+ty_0}} }{(1+ty_0)^2} \\
=&\frac{C+y_0^2\log({1+ty_0}) }{(1+ty_0)^2} \rightarrow C=0 \\ 
\frac{\partial y}{\partial \alpha_1} =& \frac{y_0^2\log({1+ty_0}) }{(1+ty_0)^2} \numberthis
\end{align*}    

\subsection{Partial derivative of higher-order $\alpha$'s}

Using similar arguments, we arrive at the conclusion that for $\alpha_n$ where $n>1$
\begin{equation}
\frac{\partial y}{\partial \alpha_n}=\frac{y_0^{n+1}((1+t y_0)^{1-n}-1)} {(1-n)(1+ty_0)^2}    
\end{equation}

Plots of the sensitivities suggest that $r$ is the dominant parameter for values of $y_0<1$, though exactly where this transition occurs is probably worth investigating.

The top-left entry in the FIM is
\begin{eqnarray}
    \begin{aligned}
     \mathcal{I}_{1,1} &= \left(\frac{\partial y}{\partial r}\right)^2 \\
 &= \left(\frac{y_0t(y_0 t+2)} {2(y_0 t+1)^2}\right)^2  \\
 &= \frac{y_0^2 t^2  (y_0 t+1)^2}{4 (y_0 t+1)^4}             
    \end{aligned}    
\end{eqnarray}

In the limit $t\rightarrow \infty$, this approaches $\frac {t^4}{t^4}$ which is order $\mathcal{O}(t^0)$. This implies that the leading eigenvector of transcritical bifurcations will be relevant, not hyperrelevant like for all other forms of bifurcations considered here. It is tempting to speculate that the topological interpretation of this quirk in the algebra stems from the unique flow-field around transcritical bifurcations. For $r<0$, the vector field has a negative-positive-negative pattern; for $r>0$ this negative-positive-negative pattern is duplicated, just with an unstable equilibrium at $y=0$ which had been stable before. Only at the critical value itself ($r=0$) is there a topological inhomogeneity. The other bifurcations have fundamentally different flow-fields on either side of the critical value, and thus, perhaps, their bifurcation parameters acquire hyper-relevance rather than simply relevance. Further study is needed to prove this conjecture.

Because $\frac{\partial y}{\partial \alpha_1} \rightarrow \mathcal{O}(\log(t)-2)$ and $\frac{\partial y}{\partial \alpha_n} \rightarrow \mathcal{O}(t^{-1-n})$, simple multiplication shows that all the other entries in the FIM will be of lower order than the top-left.

\section{\label{ap:pitch}FIM of Pitchfork Bifurcations}

In the supercritical case, the normal form is 
\begin{equation}
\frac{dy}{dt} = ry(t)-y(t)^3 +\alpha_1 y(t)^4 +\alpha_2 y(t)^5 + ...    
\end{equation}
and the subcritical case is the same except the sign on the cubic term changes. At the critical value of $\theta_i=0$, the system reduces to:

\begin{eqnarray}
\frac{dy}{dt} = -y^3 \rightarrow -\frac{dy}{y^3} = dt \rightarrow \frac{1}{2y^2}\Big|_{y_0}^{y(t)}=t\Big|_{0}^{t} \nonumber\\ 
\frac{1}{y(t)^2}-\frac{1}{y_0^2}=2t \rightarrow  \frac{1}{y(t)^2}=2t+\frac{1}{y_0^2} \nonumber \\
\rightarrow y(t)=\frac{y_0}{\sqrt{1+2t y_0^2}}    
\end{eqnarray}

Following the same logic, the formula for the subcritical case is 
\begin{equation}
    y(t)=\frac{y_0}{\sqrt{1-2t y_0^2}}
\end{equation} 
Note that this creates a potentially-problematic singularity at $t=\frac{1}{2y_0^2}$.

\subsection{Partial derivative of $r$}
Let the $\alpha_i$'s=0. The derivative of the normal form w.r.t. $r$ becomes:

\begin{eqnarray}
\frac{\partial}{\partial r}\left(\frac{\partial y}{\partial t} = ry-y^3 \right)\nonumber \\
\frac{\partial ^2 y}{\partial r \partial t}=\cancel{r\frac{\partial y}{\partial r}}+y-3y^2\frac{\partial y}{\partial r}\nonumber \\
\frac{\partial w}{\partial t}=y-3y^2w    
\end{eqnarray}

where $w=\frac{\partial y}{\partial r}$. Using integration factors $p_1=1,p_0=3y^2,q=y$, we see that 
\begin{equation*}
\begin{aligned}
\mu &= \exp\left(\int-3y^2 dt \right) \\
&= \exp\left(-\int \frac {3y_0^2 dt}{1+2y_0^2t} \right) \\
&= \exp \left(\frac{3}{2}\ln(1+2 y_0^2 t)\right) \\
&= (1+2 y_0^2 t)^{3/2}
\end{aligned}    
\end{equation*}

Therefore,
\begin{align*}
w &= \frac{C+\int{\mu y(t) dt}}{\mu}\\
&=\frac{\cancel{C}+\int\frac{y_0}{\sqrt{1+2t y_0^2}}(1+2 y_0^2 t)^{3/2} dt}{(1+2 y_0^2 t)^{3/2}} \\
\frac{\partial y}{\partial r} &=\frac{y_0 t (1+y_0^2 t)} {(1+2 y_0^2 t)^{3/2}} \numberthis\\
\end{align*}    

Following the same logic for the subcritical case eventually brings us to 
\begin{equation}
\frac{\partial y}{\partial r} = \frac{t y_0 \left(1-t y_0^2\right)}{\left(1-2 t y_0^2\right)^{3/2}}    
\end{equation}

\subsection{Partial derivative of $\alpha$'s}

When $r=0$, and all $\alpha_{i\neq n}=0$, then the normal form reduces to 
\begin{equation}
    \frac{dy}{dt} = -y(t)^3 +\alpha_n y(t)^{n+3}
\end{equation}
which conveniently allows us to use the same $\mu$ integration factor as above. Using the integration scheme outlined there, after many steps we reach the conclusion that 
\begin{equation}
    \frac{\partial y}{\partial \alpha_n} = \frac{y_0^{n+1}}{2 - n} \frac{(1 + 2 t y_0^2)^{1 - n/2} - 1}{\mu}
\end{equation}

This produces an obvious problem when $n=2$, but in that case the integration step simplifies and we find that 
\begin{equation}
\frac{\partial y}{\partial \alpha_2} = \frac{y_0^3\ln(1 + 2 t y_0^2)}{2\mu}    
\end{equation}

All this indicates that in the FIM, the entry corresponding to $(\partial y/\partial r)^2$ is $\mathcal{O}(t^1)$, while all other entries are lower order, so $r$ will be the only hyperrelevant direction. 

\section{\label{ap:hopf}FIM of Hopf Bifurcations}

Analysis of the Hopf bifurcation in either the complex or Cartesian formulation is complicated, because the introduction of nuisance parameters to the normal form equations tends to alter the period of limit cycles. This means standard trigonometric functions would also need to be altered with time-dependent terms to dilate/expand the period for a closed form solution of the trajectories $z(t)$ or $x(t), y(t)$ respectively. 

However, reparameterizing the equation into polar coordinate form simplifies matters greatly. The system $\dot{r} = r(\mu-r^2); ~\dot\theta = -1$ should look familiar, as the equation for $r$ is simply the normal form for a supercritical pitchfork bifurcation. Therefore, deriving the elements of its Fisher information matrix has already been performed in Appendix~\ref{ap:pitch}, albeit with different variable and parameter names.

\bibliography{bifurc}% Produces the bibliography via BibTeX.

\end{document}